\documentclass[letterpaper, 10 pt, conference]{ieeeconf}%
\usepackage{graphics}
\usepackage{epsfig}
\usepackage{cite}
\usepackage{breakurl}
\usepackage{amsmath}
\usepackage{amssymb}
\usepackage{amsfonts}
\usepackage{graphicx}%
\providecommand{\U}[1]{\protect\rule{.1in}{.1in}}
\IEEEoverridecommandlockouts
\newtheorem{theorem}{Theorem}[section]
\newtheorem{lemma}[theorem]{Lemma}
\newtheorem{definition}[theorem]{Definition}

\newtheorem{remark}[theorem]{Remark}

\newtheorem{assumption}[theorem]{Assumption}

\begin{document}

\title{{\LARGE \textbf{Optimal Control and Coordination of Connected and Automated
Vehicles at Urban Traffic Intersections }}}
\author{Yue J. Zhang, Andreas A. Malikopoulos, Christos G. Cassandras \thanks{This
manuscript has been authored by UT-Battelle, LLC under Contract No.
DE-AC05-00OR22725 with the U.S. Department of Energy. The United States
Government retains and the publisher, by accepting the article for
publication, acknowledges that the United States Government retains a
non-exclusive, paid-up, irrevocable, world-wide license to publish or
reproduce the published form of this manuscript, or allow others to do so, for
United States Government purposes.} \thanks{This research was supported by the
Laboratory Directed Research and Development Program of Oak Ridge National
Laboratory, managed by UT-Battelle, LLC, for the U. S. Department of Energy.
The work of Cassandras and Zhang is supported in part by NSF under grants CNS-
1239021, ECCS-1509084, and IIP-1430145, by AFOSR under grant FA9550-15-1-0471,
and by ONR under grant N00014-09-1-1051.} \thanks{Y.J. Zhang and C.G.
Cassandras are with the Division of Systems Engineering and Center for
Information and Systems Engineering, Boston University, Boston, MA 02215 USA
(e-mail: joycez@bu.edu; cgc@bu.edu).} \thanks{A.A. Malikopoulos is with the
Energy \& Transportation Science Division, Oak Ridge National Laboratory, Oak
Ridge, TN 37831 USA (e-mail: andreas@ornl.gov).} }
\maketitle

\begin{abstract}
We address the problem of coordinating online a continuous flow of connected
and automated vehicles (CAVs) crossing two adjacent intersections in an urban
area. We present a decentralized optimal control framework whose solution
yields for each vehicle the optimal acceleration/deceleration at any time in
the sense of minimizing fuel consumption. The solution, when it exists, allows
the vehicles to cross the intersections without the use of traffic lights,
without creating congestion on the connecting road, and under the hard safety
constraint of collision avoidance. The effectiveness of the proposed solution
is validated through simulation considering two intersections located in
downtown Boston, and it is shown that coordination of CAVs can reduce
significantly both fuel consumption and travel time.

\end{abstract}

\thispagestyle{empty} \pagestyle{empty}

%%%%%%%%%%%%%%%%%%%%%%%%%%%%%%%%%%%%%%%%%%%%%%%%%%%%%%%%%%%%%%%%%%%%%%%%%%%%%%%%

%%%%%%%%%%%%%%%%%%%%%%%%%%%%%%%%%%%%
%SECTION I: INTRODUCTION
%%%%%%%%%%%%%%%%%%%%%%%%%%%%%%%%%%%%

\section{Introduction}

Connected and automated vehicles (CAVs) can improve transportation safety and
efficiency using traffic lights and vehicle-to-infrastructure communication
\cite{Li2014}. There are also significant opportunities  to coordinate CAVs
for improving both safety and traffic flow using either centralized or
decentralized approaches. In this paper, we categorize an approach as
centralized if there is at least one task in the system that is globally
decided for all vehicles by a single central controller. In a decentralized
approach, a \textquotedblleft coordinator\textquotedblright\ may be used to
handle or distribute information available in the system without, however,
getting involved in any control task.

To date, traffic lights are the prevailing method used to control the traffic
flow through an intersection. Recent technological developments which exploit
the ability to collect traffic data in real time have made it possible for new
methods to be applied to traffic light control \cite{Liu2007}. Most of these
approaches are computationally inefficient and not immediately amenable to
online implementations. More recently, however, data-driven approaches have
been developed leading to on-line adaptive traffic light control as in
\cite{Fleck2015}. Aside from the obvious infrastructure cost and the need for
dynamically controlling green/red cycles, traffic light systems also lead to
problems such as significantly increasing the number of rear-end collisions at
an intersection. These issues have provided the motivation for drastically new
approaches capable of providing a smoother traffic flow and more
fuel-efficient driving while also improving safety.

The advent of CAVs provides the opportunity for such new approaches. Dresner
and Stone \cite{Dresner2004} proposed a reservation scheme for automated
vehicle intersection control whereby a centralized controller coordinates a
crossing schedule based on requests and information received from the vehicles
located inside some communication range. This scheme has been expanded since
then \cite{Dresner2008,DeLaFortelle2010,Huang2012}. Increasing the throughput
of an intersection is one desired goal and it can be achieved through the
travel time optimization for all vehicles located within a radius from the
intersection. There have been several research efforts to address the problem
of vehicle coordination at intersections within a decentralized control
framework \cite{Milanes2010,Alonso2011,Kim2014,Campos2014}. One of the main
challenges in this case is the possibility of having deadlocks in the
solutions as a consequence of the use of local information.

In this paper, we address the problem of optimally controlling online the fuel
consumption of a varying number of CAVs subject to congestion and safety
constraints as they cross two urban intersections. The contribution of the
paper is a decentralized control problem framework whose solution yields for
each vehicle the optimal acceleration/deceleration at any time without
creating congestion on the connecting road and under the hard constraint of
collision avoidance.

The structure of the paper is as follows. In Section II, we extend our work on
a single intersection \cite{Malikopoulos2016} and provide a model for two
intersections. In Section III, we formulate the problem of CAV coordination
and optimal control for two intersections and provide an analytical solution.
In Section IV, we present simulation results in theVISSIM simulation
environment considering two intersections located in downtown Boston and offer
concluding remarks in Section V.

%%%%%%%%%%%%%%%%%%%%%%%%%%%%%%%%%%%%
%SECTION II: The Model
%%%%%%%%%%%%%%%%%%%%%%%%%%%%%%%%%%%%

\section{The Model}

\begin{figure}[ptb]
\centering
\includegraphics[width=3.2 in]{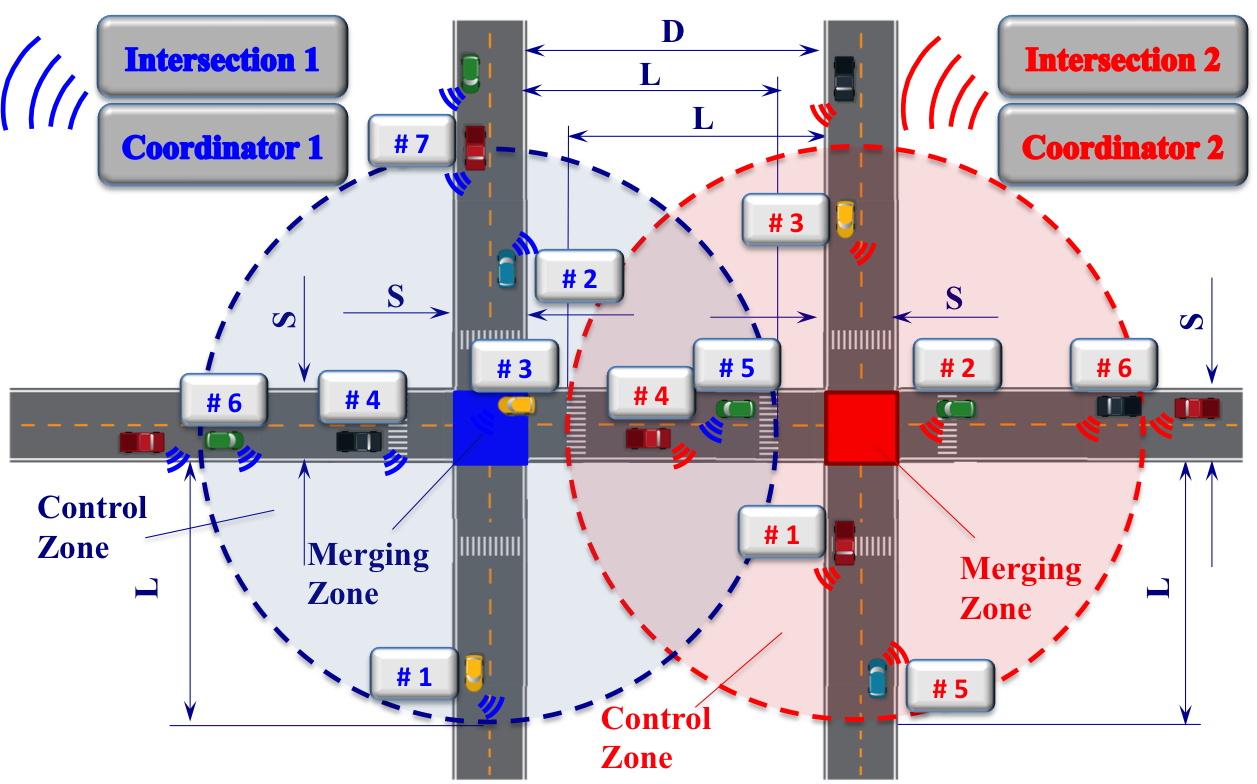} \caption{Two intersections
with connected and automated vehicles.}%
\label{fig:intersection}%
\end{figure}

We consider two intersections, 1 and 2, located within a distance $D$ (Fig.
\ref{fig:intersection}). The region at the center of each intersection, called
\textit{merging zone}, is the area of potential lateral collision of the
vehicles. Although this is not restrictive, we consider the \textit{merging
zones} in both intersections to be squares of equal sides $S$. Each
intersection has a \textit{control zone} and a coordinator that can
communicate with the vehicles traveling within it. The distance between the
entry of the control zone and the entry of the merging zone is $L>S$, and it
is assumed to be the same for all entry points to a given control zone.

We consider a time-varying number of CAVs $N_{z}(t)\in\mathbb{N}$ present at
control zone $z=1,2$ at time $t\in\mathbb{R}$. When a CAV reaches the control
zone of intersection $z$ at some instant $t$, the coordinator assigns a unique
identity consisting of a pair $(i,j)$. Here, $i=N_{z}(t)+1$ is an integer
corresponding to the position of the CAV in a first-in-first-out (FIFO) queue
for this control zone. The elements of this queue can belong to any of four
subsets (precisely defined in Definition \ref{def:2}) depending on the road
and lane traveled by each CAV so that $j\in\{1,\ldots,4\}$ is an integer
corresponding to the appropriate subset. If two or more vehicles enter the
control zone of any intersection at the same time, then the corresponding
coordinator selects randomly their position in the queue.

The vehicles in the control zone of intersection $z=1$ traveling from west to
east (see Fig. \ref{fig:intersection}) remain in the queue imposed by
coordinator 1 until they exit the corresponding merging zone. In the region
between the exit point of merging zone 1 and the entry point of control zone
2, the vehicles cruise with the speed they had when they exited that merging
zone and then enter the queue imposed by the coordinator of intersection
$z=2$. A similar process applies to vehicles in control zone 2 traveling from
east to west.

The objective of each vehicle is to derive an optimal
acceleration/deceleration at any time so as to minimize fuel consumption over
the time interval defined from its entry time at a control zone to its exit
time from the merging zone while avoiding congestion between the two
intersections. We consider an indication of potential congestion the speed
reduction of any of the vehicles traveling on this road below a desired
minimum value. Accordingly, we specify congestion-avoidance constraints as
described in the next section. In addition, we impose hard constraints so to
avoid either rear-end collision, or lateral collision inside the merging zone.

Let $\mathcal{N}_{z}(t)=\{1,\ldots,N_{z}(t)\},$ $z=1,2$, be the queue
associated with the control zone of intersection $z$. We represent the
dynamics of each vehicle $i$, $i\in\mathcal{N}_{z}(t)$, moving along a
specified lane with a state equation
\begin{equation}
\dot{x_{i}}=f(t,x_{i},u_{i}),\qquad x_{i}(t_{i}^{0})=x_{i}^{0}%
,\label{eq:model}\\
\end{equation}
where $t\in\mathbb{R}^{+}$ is the time, $x_{i}(t)$, $u_{i}(t)$ are the state
of the vehicle and control input, $t_{i}^{0}$ is the time that vehicle $i$
enters the control zone, and $x_{i}^{0}$ is the value of the initial state.
For simplicity, we assume that each vehicle is governed by a second order
dynamics
\begin{equation}
\dot{p}_{i}=v_{i}(t) ~ \text{,} ~ \dot{v}_{i}=u_{i}(t) \label{eq:model2}%
\end{equation}
where $p_{i}(t)\in\mathcal{P}_{i}$, $v_{i}(t)\in\mathcal{V}_{i}$, and
$u_{i}(t)\in\mathcal{U}_{i}$ denote the position, speed and
acceleration/deceleration (control input) of each vehicle $i$. Let
$x_{i}(t)=\left[
\begin{array}
[c]{cc}%
p_{i}(t) & v_{i}(t)
\end{array}
\right]  ^{T}$ denote the state of each vehicle $i$, with initial value
$x_{i}^{0}=\left[
\begin{array}
[c]{cc}%
0 & v_{i}^{0}%
\end{array}
\right]  ^{T}$, taking values in the state space $\mathcal{X}_{i}%
=\mathcal{P}_{i}\times\mathcal{V}_{i}$. The sets $\mathcal{P}_{i}$,
$\mathcal{V}_{i}$ and $\mathcal{U}_{i}$, $i\in\mathcal{N}_{z}(t),$ $z=1,2,$,
are complete and totally bounded subsets of $\mathbb{R}$. It follows that
$\mathcal{P}_{i}$, $\mathcal{V}_{i}$, and $\mathcal{U}_{i}$ are Borel
measurable sets. The state space $\mathcal{X}_{i}$ for each vehicle $i$ is
closed with respect to the induced topology on $\mathcal{P}_{i}\times
\mathcal{V}_{i}$ and thus, it is compact.

We need to ensure that for any initial state $(t_{i}^{0},x_{i}^{0})$ and every
admissible control $u(t)$, the system \eqref{eq:model} has a unique solution
$x(t)$ on some interval $[t_{i}^{0},t_{i}^{f}]$, where $t_{i}^{f}$ is the time
that vehicle $i\in\mathcal{N}_{z}(t)$ exits the merging zone of intersection
$z$. The following observations from \eqref{eq:model2} satisfy some regularity
conditions required both on $f$ and admissible controls $u(t)$ to guarantee
local existence and uniqueness of solutions for \eqref{eq:model}:
%%%%%%%%
%Observations
%%%%%%%%
a) $f$ is continuous in $u$ and continuously differentiable in $x$, b) The
first derivative of $f$ in $x$, $f_{x}$, is continuous in $u$, and c) The
admissible control $u(t)$ is continuous in $t$. We impose the following
assumption regarding the final conditions when a vehicle exits the merging
zone,
%%%%%%%%
%Assumption
%%%%%%%%
which is intended to enhance safety awareness:

\begin{assumption}
\label{ass:4} The vehicle speed inside any merging zone is constant.
\end{assumption}

This assumption is not restrictive and could be modified appropriately. In
addition, to ensure that the control input and vehicle speed are within a
given admissible range, the following constraints are imposed.
\begin{equation}%
\begin{split}
u_{min} &  \leq u_{i}(t)\leq u_{max},\quad\text{and}\\
0 &  \leq v_{min}\leq v_{i}(t)\leq v_{max},\quad\forall t\in\lbrack t_{i}%
^{0},t_{i}^{f}],
\end{split}
\label{speed_accel constraints}%
\end{equation}
where $u_{min}$, $u_{max}$ are the minimum deceleration and maximum
acceleration allowable, and $v_{min}$, $v_{max}$ are the minimum and maximum
speed limits respectively.

%%%%%%%%%%%%%%%%%%%%%%%%%%%%%%%%%%%%%%%%%%%%%%%%%%%%%%%%%%%%%%%%%%%%%%%%%%%%%%%%
%%%%%%%%%%%%%%%%%%%%%%%%%%%%%%%%%%%%
%SECTION III: Vehicle Coordination
%%%%%%%%%%%%%%%%%%%%%%%%%%%%%%%%%%%%

\section{Vehicle Coordination}

%%%%%%%%%%%%%%%%%%%%%%%%%%%%%%%%%%%%%%%
%Decentralized Problem Formulation
%%%%%%%%%%%%%%%%%%%%%%%%%%%%%%%%%%%%%%%

\subsection{Decentralized Control Problem Formulation}

When a vehicle enters a control zone $z=1,2$, it receives a unique identity
$(i,j)$ from the coordinator, as described in the previous section. Since the
coordinator is not involved in any decision on the vehicle control, we can
formulate $N_{1}(t)$ and $N_{2}(t)$ decentralized tractable problems for
intersection 1 and 2 respectively that may be solved on line. Before we
proceed with the decentralized problem formulation we need to establish some definitions.

Recall that $\mathcal{N}_{z}(t)=\{1,\ldots,N_{z}(t)\}$ is the FIFO queue of
vehicles in control zone $z=1,2$. A vehicle index $i\in\mathcal{N}_{z}(t)$
also indicates which vehicle is closer to the merging zone, i.e., if $i<k$
then $L-p_{i}<L-p_{k}$.

%%%%%%%%
%Definition
%%%%%%%%

\begin{definition}
\label{def:2} Each vehicle $i\in\mathcal{N}_{z}(t)$ belongs to at least one of
the following four subsets: 1) $\mathcal{R}_{i}^{z}(t)$ contains all vehicles
traveling on the same road as vehicle $i$ and towards the same direction but
on different lanes, 2) $\mathcal{L}_{i}^{z}(t)$ contains all vehicles
traveling on the same road and lane as vehicle $i$, 3) $\mathcal{C}_{i}%
^{z}(t)$ contains all vehicles traveling on different roads from $i$ and
having destinations that can cause collision at the merging zone, and 4)
$\mathcal{O}_{i}^{z}(t)$ contains all vehicles traveling on the same road as
vehicle $i$ and opposite destinations that cannot, however, cause collision at
the merging zone.
\end{definition}

%%%%%%%%
%Definition
%%%%%%%%

To illustrate the definitions of the subsets of $\mathcal{N}_{z}(t)$, observe
that in Fig. \ref{fig:intersection} vehicles \# 4 and \# 6 (blue label) belong
to $\mathcal{L}_{6}^{1}(t)$; vehicles \# 4 and \# 7 (blue label) belong to
$\mathcal{C}_{7}^{1}(t)$ while vehicles \# 2 and \# 3 (red label) belong to
$\mathcal{C}_{3}^{2}(t)$; vehicles \# 4 and \# 5 (blue label) belong to
$\mathcal{O}_{5}^{1}(t)$ while vehicles \# 3 and \# 5 (red label) belong to
$\mathcal{O}_{5}^{2}(t)$.

\begin{definition}
\label{def:id} The \textit{unique identity} that the coordinator assigns to
each vehicle $i\in\mathcal{N}_{z},$ $z=1,2,$ at time $t$ when the vehicle
arrives at control zone $z$, is a pair $(i,j)$ where $i=N_{z}(t)+1$ is an
integer representing the location of the vehicle in the FIFO queue
$\mathcal{N}_{z}(t)$ and $j\in\{1,\ldots,4\}$ is an integer based on a
one-to-one mapping from $\{\mathcal{R}_{i}^{z}(t),$ $\mathcal{L}_{i}^{z}(t),$
$\mathcal{C}_{i}^{z}(t),$ $\mathcal{O}_{i}^{z}(t)\}$ onto $\{1,\ldots,4\}$.
\end{definition}

%%%%%%%%
%Assumption
%%%%%%%%

\begin{assumption}
\label{ass:sensor} Each vehicle $i$ has proximity sensors and can observe
and/or estimate local information that can be shared with other vehicles.
\end{assumption}

%%%%%%%%
%Definition
%%%%%%%%

\begin{definition}
\label{def:4} For each vehicle $i$ when it enters a control zone, we define
the \textit{information set} $Y_{i}(t)$ as
\begin{equation}%
\begin{split}
Y_{i}(t) &  \triangleq\Big\{p_{i}(t),v_{i}(t),\mathcal{Q}_{j}^{z}%
,j=1,\ldots,4,z=1,2,s_{i}(t),t_{i}^{f}\Big\},\\
\forall t &  \in\lbrack t_{i}^{0},t_{i}^{f}],
\end{split}
\end{equation}

\end{definition}
where $p_{i}(t),v_{i}(t)$ are the position and speed of vehicle $i$ inside the
\textit{control zone }it belongs to, and $\mathcal{Q}_{j}^{z}\in
\{\mathcal{R}_{i}^{z}(t),$ $\mathcal{L}_{i}^{z}(t),$ $\mathcal{C}_{i}^{z}(t),$
$\mathcal{O}_{i}^{z}(t)\},$ $z=1,2,$ is the subset assigned to vehicle $i$ by
the coordinator (see Definition \ref{def:2}). The first of the two new
elements in $Y_{i}(t)$ yet to be defined is $s_{i}(t)=p_{k}(t)-p_{i}(t)$; this
represents the distance between vehicle $i$ and some vehicle $k$ which is
immediately ahead of $i$ in the same lane (the index $k$ is made available to
$i$ by the coordinator). The last element above, $t_{i}^{f}$, is the time
targeted for vehicle $i$ to exit the merging zone, whose evaluation is
discussed next. Note that once the vehicle $i$ enters the control zone, then
immediately all information in $Y_{i}(t)$ becomes available to $i$:
$p_{i}(t),v_{i}(t)$ are read from the sensors; $\mathcal{Q}_{j}^{z}$ is
assigned by the coordinator, as is the value of $k$ based on which $s_{i}(t)$
is also evaluated; $t_{i}^{f}$ can also be computed at that time, as described next.

The time $t_{i}^{f}$ that the vehicle $i$ exits the merging zone is based on
imposing constraints aimed at avoiding congestion (in the sense of maintaining
vehicle speeds above a certain value). There are three cases to consider,
depending on the value of $\mathcal{Q}_{j}^{z}$:

1) if the predecessor of vehicle $i$ in queue $\mathcal{N}_{z}(t)$, i.e.,
vehicle $i-1$, belongs to either $\mathcal{R}_{i}^{z}(t)$ or $\mathcal{O}%
_{i}^{z}(t),$ $z=1,2,$ then both $i-1$ and $i$ can share the merging zone at
the same time; thus, to minimize the distances between vehicles in the queue
(hence, not unnecessarily reduce speeds) both $i-1$ and $i$ should be entering
and exiting the merging zone at the same time. Therefore, we impose the
constraint $t_{i}^{f}=t_{i-1}^{f}$.

2) If vehicle $i-1$ belongs to $\mathcal{L}_{i}^{z}(t),$ $z=1,2,$ then, by the
same argument, both $i-1$ and $i$ should have the \emph{minimal safe distance}
allowable, denoted by $\delta$, by the time vehicle $i-1$ enters the merging
zone, i.e., $t_{i}^{f}=t_{i-1}^{f}+\frac{\delta}{v_{i}(t_{i-1}^{f}%
)}\label{tif_2},$ where $v_{i}(t_{i-1}^{f})=v_{i-1}(t_{i-1}^{f})$.

3) Finally, if vehicle $i-1$ belongs to $\mathcal{C}_{i}^{z}(t),$ $z=1,2,$ we
constrain the merging zone to contain only one vehicle so as to avoid a
lateral collision. Therefore, vehicle $i$ is allowed to enter the merging zone
only when vehicle $i-1$ exits the merging zone, where $t_{i}^{m}$ is the time
that the vehicle $i$ enters the merging zone), i.e., $t_{i}^{f}=t_{i-1}%
^{f}+\frac{S}{v_{i}(t_{i-1}^{f})}\label{tif_3},$ where $v_{i}(t_{i-1}%
^{f})=\frac{L}{t_{i-1}^{f}-t_{i}^{0}}$.

Note that, in all cases, once vehicle $i$ enters the control zone, vehicle
$i-1$ is already present, thus $t_{i-1}^{f}$, $v_{i-1}(t_{i-1}^{f})$, and $\mathcal{Q}_j^z$, $z=1,2,$ are available through $Y_{i-1}(t)$. Moreover, to
ensure the absence of rear-end collision between two consecutive vehicles
traveling on the same lane we impose the constraint $s_{i}(t)\geq\delta$
(obviously, this applies only when $N_{z}(t)>1$).

However, $t_{i}^{f}$ above may not be feasible due to the speed and
acceleration constraint in \eqref{speed_accel constraints}. There are two
cases to consider, based on whether vehicle $i$ can reach $v_{max}$ prior to
$t_{i-1}^{f}$: 

$(i)$ If vehicle $i$ enters the control zone at $t_{i}^{0}$, it
accelerates with $u_{max}$ until it reaches $v_{max}$ and then cruises at this
speed until it leaves the merging zone at time $t_{i}^{1}$. It is
straightforward to show (details found in
\cite{ZhangMalikopoulosCassandras2016}) that
\begin{equation}
t_{i}^{1}=t_{i}^{0}+\frac
{L+S}{v_{max}}+\frac{(v_{max}-v_{i}^{0})^{2}}{2u_{max}v_{max}}.
\end{equation}
 
$(ii)$ Vehicle $i$ reaches the merging zone at $t_{i}^{m}$ with speed $v_{i}%
(t_{i}^{m})<v_{max}$. It is again straightforward to show that in this case
\begin{equation}
t_{i}^{2}=t_{i}^{0}+\frac{v_{i}(t_{i}^{m})-v_{i}^{0}}{u_{max}}+\frac
{S}{v_{i}(t_{i}^{m})}.
\end{equation} 
where $v_{i}(t_{i}^{m})=\sqrt{2Lu_{max}+(v_{i}%
^{0})^{2}}$. Thus, setting $t_{i}^{c}=\max\{t_{i}^{1},t_{i}^{2}\}$, the value
of $t_{i}^{f}$ is computed as described above and summarized as follows, where
$z=1,2$:%
\begin{equation}
t_{i}^{f}=\left\{
\begin{array}
[c]{ll}%
t_{1}^{f}, & \mbox{if $i=1$},\\
\text{max }\{t_{i-1}^{f},t_{i}^{c}\}, &
\mbox{if $i-1\in\mathcal{R}_{i}^{z}$ (or $\mathcal{O}_{i}^{z}),$}\\
\text{max }\{t_{i-1}^{f}+\frac{\delta}{v_{i}(t_{i-1}^{f})},t_{i}^{c}\}, &
\mbox{if $i-1\in\mathcal{L}_{i}^{z}$},\\
\text{max }\{t_{i-1}^{f}+\frac{S}{v_{i}(t_{i-1}^{f})},t_{i}^{c}\}, &
\mbox{if $i-1\in\mathcal{C}_{i}^{z}$}.
\end{array}
\right.  \label{def:tf}%
\end{equation}

%%%%%%%%
%Definition
%%%%%%%%

\begin{definition}
\label{def:rear} For each vehicle $i\in\mathcal{N}_{z}(t), z=1,2,$ we define
the \textit{rear-end control interval} $R_{i}$ as
\begin{gather}
R_{i}(t) \triangleq\Big\{u_{i}(t) \in[u_{min}, u_{max}] ~|~ s_{i}(t) \ge
\delta,\nonumber\\
\forall i\in\mathcal{N}_{z}(t), z=1,2, \forall t \in\lbrack t_{i}^{0}%
,t_{i}^{f}], |\mathcal{N}_{z}|>1\Big\}.\label{eq:moddef1}%
\end{gather}

\end{definition}

%%%%%%%%
%Remark
%%%%%%%%

\begin{remark}
\label{def:remark} At each time $t$, each vehicle $i\in\mathcal{N}_{z},z=1,2,$
communicates with the preceding vehicle $i-1$ in the queue and accesses the
values of $t_{i-1}^{f}$, $v_{i-1}(t_{i-1}^{f})$, $\mathcal{Q}_j^z$,
$j=1,\dots,4$, $z=1,2$ \textit{from its information set} (Definition \ref{def:4}).
This information is necessary for vehicle $i$ to compute $t_{i}^{f}$
appropriately and satisfy \eqref{def:tf} and \eqref{eq:moddef1}.
\end{remark}

%%%%%%%%
%Lemma 1
%%%%%%%%

\begin{lemma}
The decentralized communication structure aims for each vehicle $i$ to solve
an optimal problem for $t\in[t^{0}_{i}, t^{f}_{i}]$ the solution of which
depends only on the solution of the vehicle $i$-1.
\end{lemma}

Due to space limitations, all proofs are omitted but may be found in
\cite{ZhangMalikopoulosCassandras2016}.

Consequently the decentralized control problem for each CAV approaching either
intersection can be formulated so as to minimize the $L^{2}$-norm of its
control input (acceleration/deceleration). It has been shown
\cite{Rios-Torres2015} that there is monotonic relationship between fuel
consumption for each vehicle $i$, $f^{fuel}_{i}(t)$, and its control input
$u_{i}$. Thus, the problem of minimizing the acceleration/deceleration is
equivalent to the problem of minimizing fuel consumption, and it is formulated
as follows:
\begin{gather}
\min_{u_{i}\in R_{i}} \frac{1}{2} \int_{t^{0}_{i}}^{t^{f}_{i}} u^{2}_{i} \cdot
K_{i} ~ dt\nonumber\\
\text{Subject to}: \eqref{eq:model2}, \eqref{def:tf} ~ \forall i
\in\mathcal{N}_{z}, z= 1, 2, \label{eq:decentral}%
\end{gather}
where $K_{i}$ is a factor to capture CAV diversity. However, for simplicity in
the rest of the paper we set $K_{i}=1$. Both rear-end and lateral collision
avoidance constraints are satisfied at time $t^{f}_{i}$.

%%%%%%%%%%%%%%%%%%%%%%%%%%%%%%%%%%%%%%%
%Analytical solution of the decentralized control problem
%%%%%%%%%%%%%%%%%%%%%%%%%%%%%%%%%%%%%%%

\subsection{Analytical solution of the decentralized control problem}

For the analytical solution and online implementation of the decentralized
problem \eqref{eq:decentral}, we apply Hamiltonian analysis by considering
that when the CAVs enter the control zone, the constraints are not active.
Clearly, this is in general not true. For example, a vehicle may enter the
control zone with speed higher than the speed limit. In this case, we need to
solve an optimal control problem starting from an infeasible state. To address
this situation requires additional analysis which is the subject of ongoing research.

From \eqref{eq:decentral}, the state equations \eqref{eq:model2}, the
control/state constraints \eqref{speed_accel constraints}, and rear-end
collision avoidance constraint for each vehicle $i\in\mathcal{N}%
_{z}(t),z=1,2,$ the Hamiltonian function can be formulated as follows
\begin{gather}
H_{i}\big(t,x(t),u(t)\big)=\frac{1}{2}u_{i}^{2}+\lambda_{i}^{p}\cdot
v_{i}+\lambda_{i}^{v}\cdot u_{i}\nonumber\\
+\mu_{i}^{a}\cdot(u_{i}-u_{max})+\mu_{i}^{b}\cdot(u_{min}-u_{i})+\mu_{i}%
^{c}\cdot(v_{i}-v_{max})\nonumber\\
+\mu_{i}^{d}\cdot(v_{min}-v_{i})+\mu_{i}^{\delta}\cdot(p_{i}-p_{k}%
+\delta),i,k\in\mathcal{N}_{z},z=1,2,\nonumber\label{eq:16}%
\end{gather}
where $\lambda_{i}^{p}$ and $\lambda_{i}^{v}$ are the co-state components, and
$\mu^{T}$ is the vector of the Lagrange multipliers. The solution of the
problem including the rear-end collision avoidance constraint may become
intractable due to the numerous scenarios of activation/deactivation of the
constraints. Thus, we will not include it in the analysis below. Note that we
can guarantee rear-end collision avoidance at time $t_{i}^{f}$, but it remains
to show that the constraint does not become active at any time in $(t_{i}%
^{0},t_{i}^{f}]$ assuming it is not active at $t=t_{i}^{0}$. The Lagrange
multipliers are $\mu_{i}^{a}=\mu_{i}^{b}=\mu_{i}^{c}=\mu_{i}^{d}=0$ if the
constraints are not active, and they are greater than zero if the constraints
become active. The Euler-Lagrange equations become
\begin{equation}
\dot{\lambda}_{i}^{p}=-\frac{\partial H_{i}}{\partial p_{i}}=0,\label{eq:EL1}%
\end{equation}
and
\begin{equation}
\dot{\lambda}_{i}^{v}=-\frac{\partial H_{i}}{\partial v_{i}}=\left\{
\begin{array}
[c]{ll}%
-\lambda_{i}^{p}, & \mbox{$v_{i}(t) - v_{max} <0$}~\text{and}\\
& \mbox{$v_{min} - v_{i}(t)>0$},\\
-\lambda_{i}^{p}+\mu_{i}^{c}, & \mbox{$v_{i}(t) - v_{max} =0$},\\
-\lambda_{i}^{p}-\mu_{i}^{d}, & \mbox{$v_{min} - v_{i}(t)=0$}.
\end{array}
\right.  \label{eq:EL2}%
\end{equation}
The necessary condition for optimality is%

\begin{equation}
\label{eq:KKT1}\frac{\partial H_{i}}{\partial u_{i}} = u_{i} + \lambda^{v}_{i}
+ \mu^{a}_{i} - \mu^{b}_{i} = 0,
\end{equation}

To address this problem, the constrained and unconstrained arcs will be pieced
together to satisfy the Euler-Lagrange equations and necessary condition of optimality.

%%%%%%%%%%%%%%%%%%%%%%%%%%%%%%%%%%
%Case 1
%%%%%%%%%%%%%%%%%%%%%%%%%%%%%%%%%%

If the inequality control and state constraints are not active, we have
$\mu^{a}_{i} = \mu^{b}_{i}= \mu^{c}_{i}=\mu^{d}_{i}=\mu^{\delta}_{i}0.$
Applying the necessary condition \eqref{eq:KKT1}, the optimal control can be
given
\begin{equation}
\label{eq:17}u_{i} + \lambda^{v}_{i}= 0, \quad i \in\mathcal{N}(t).\\
\end{equation}
The Euler-Lagrange equations yield
%\begin{equation}
$\dot\lambda^{p}_{i} = - \frac{\partial H_{i}}{\partial p_{i}} =
0\label{eq:18}$
%\end{equation}
%\begin{equation}
and $\dot\lambda^{v}_{i} = - \frac{\partial H_{i}}{\partial v_{i}} = -
\lambda^{p}_{i};
%\label{eq:19}%
$
%\end{equation}
hence $\lambda^{p}_{i} = a_{i}$ and $\lambda^{v}_{i} = -(a_{i} t+ b_{i})$,
where $a_{i}$ and $b_{i}$ are constants of integration corresponding to each
vehicle $i$. Consequently, the optimal control input
(acceleration/deceleration) as a function of time is given by
\begin{equation}
\label{eq:20}u^{*}_{i}(t) = a_{i}t +b_{i}.\\
\end{equation}
Substituting the last equation into the vehicle dynamics equations
\eqref{eq:model2} we can find the optimal speed and position for each vehicle,
namely
\begin{equation}
v^{*}_{i}(t) = \frac{1}{2} a_{i}t^{2} + b_{i}t + c_{i} \label{eq:21}%
\end{equation}
\begin{equation}
p^{*}_{i}(t) = \frac{1}{6} a_{i}t^{3} + \frac{1}{2} b_{i}t^{2} + c_{i}t +
d_{i}, \label{eq:22}%
\end{equation}
where $c_{i}$ and $d_{i}$ are constants of integration. These constants can be
computed by using the initial and final conditions. Since we seek to derive
the optimal control (\ref{eq:20}) online, we can designate initial values
$p_{i}(t^{0}_{i})$ and $v_{i}(t^{0}_{i})$, and initial time, $t^{0}_{i}$, to
be the current values of the states $p_{i}(t)$ and $v_{i}(t)$ and time $t$,
where $t^{0}_{i} \le t \le t^{f}_{i}$. Therefore the constants of integration
will be functions of time and states, i.e., $a_{i}(t, p_{i}, v_{i}), b_{i}(t,
p_{i}, v_{i}), c_{i}(t, p_{i}, v_{i})$, and $d_{i}(t, p_{i}, v_{i})$. To
derive online the optimal control for each vehicle $i$, we need to update the
integration constants at each time $t$. Equations \eqref{eq:21} and
\eqref{eq:22}, along with the initial and final conditions defined above, can
be used to form a system of four equations of the form $\mathbf{T}_{i}
\mathbf{b}_{i} = \mathbf{q}_{i}$, namely%

\begin{equation}
\left[
\begin{array}
[c]{cccc}%
\frac{1}{6}t^{3} & \frac{1}{2}t^{2} & t & 1\\
\frac{1}{2}t^{2} & t & 1 & 0\\
\frac{1}{6}(t_{i}^{f})^{3} & \frac{1}{2}(t_{i}^{f})^{2} & t_{i}^{f} & 1\\
\frac{1}{2}(t_{i}^{f})^{2} & t_{i}^{f} & 1 & 0
\end{array}
\right]  .\left[
\begin{array}
[c]{c}%
a_{i}\\
b_{i}\\
c_{i}\\
d_{i}%
\end{array}
\right]  =\left[
\begin{array}
[c]{c}%
p_{i}(t)\\
v_{i}(t)\\
p_{i}(t_{i}^{f})\\
v_{i}(t_{i}^{f})
\end{array}
\right]  .\label{eq:23}%
\end{equation}
Hence we have
\begin{equation}
\mathbf{b}_{i}(t,p_{i}(t),v_{i}(t))=(\mathbf{T}_{i})^{-1}.\mathbf{q}%
_{i}(t,p_{i}(t),v_{i}(t)),\label{eq:24}%
\end{equation}
where $\mathbf{b}_{i}(t,p_{i}(t),v_{i}(t))$ contains the four integration
constants $a_{i}(t,p_{i},v_{i}),b_{i}(t,p_{i},v_{i}),c_{i}(t,p_{i}%
,v_{i}),d_{i}(t,p_{i},v_{i})$. Thus \eqref{eq:20} can be written as
\begin{equation}
u_{i}^{\ast}(t,p_{i}(t),v_{i}(t))=a_{i}(t,p_{i}(t),v_{i}(t))t+b_{i}%
(t,p_{i}(t),v_{i}(t)).\label{eq:25}%
\end{equation}
Since \eqref{eq:23} can be computed on line, the controller can yield the
optimal control on line for each vehicle $i$, with feedback indirectly
provided through the re-calculation of the vector $\mathbf{b}_{i}%
(t,p_{i}(t),v_{i}(t))$ in \eqref{eq:24}.

Similar results are obtained when the constraints become active. Due to space
limitations, this analysis is omitted but may be found in
\cite{Malikopoulos2016}. Note that the control for vehicle $i$ remains
unchanged until an ``event" occurs that affects its behavior. Therefore, the
time-driven controller above can be replaced by an event-driven one without
affecting its optimality properties under conditions described in
\cite{Zhong2010}.

%%%%%%%%%%%%%%%%%%%%%%%%%%%%%%%%%%%%%%%
%Interdependence of the Intersections
%%%%%%%%%%%%%%%%%%%%%%%%%%%%%%%%%%%%%%%

\subsection{Interdependence of the Intersections}

The two intersections are interdependent, i.e., the coordination of vehicles
at the merging zone of one intersection affects the behavior of vehicle
coordination of the other merging zone, and a potential congestion on the
connecting road of length $D$ (Fig. \ref{fig:intersection}) can disturb the
traffic flow. As the number of vehicles $N_{z}(t), z=1,2,$ inside the control
zones increases, the imposed safety constraints may reinforce some of the
vehicles to slow down. When the speed of a vehicle $i$, traveling on the road
that connects the two intersections, drops below the desired minimum speed,
$v_{min}$, we are interested in a control ``mechanism" to accelerate the
preceding vehicles to create more space on the road for the following vehicles.

%%%%%%%%
%Definition
%%%%%%%%

\begin{definition}
\label{def:tau} For each vehicle $i$, we define $\tau_{i}$ the additional
minimum time required for the vehicle to reach the merging zone with the
desired minimum speed $v_{min}$, i.e., $\tau_{i} = \frac{L-p_{i}(t)}{v_{min}%
}-\frac{L-p_{i}(t)}{v_{i}(t)}$.
\end{definition}

Therefore if a vehicle $i\in\mathcal{N}_{z}, z=1,2,$ travels towards an
intersection and the speed drops below the desired minimum speed
$v_{i}(t)<v_{min}$, then the first vehicle in the queue must expedite its
time, $t_{1}^{f}$, to exit the merging zone by $\tau_{i}(t)$, which means now
we have $t_{1}^{f} = t^{f}_{1}-\tau_{k}$. By doing so, it will also change
the value of $t_{2}^{f}$ of the second vehicle by $\tau_{i}(t)$ and so on, and
thus all vehicles from 1 to $i$-1 will accelerate to create the required space
for vehicle $i$ to cruise with at least $v_{min}$.

In this context, the \textit{information set} for each vehicle is expanded to
include $\tau_{i}$, namely
\begin{gather}
Y_{i}(t) \triangleq\Big\{ p_{i}(t), v_{i}(t), \mathcal{Q}_{j}^{z}, j=1,\ldots,
4, z=1,2, s_{i}(t), t_{i}^{f}, \tau_{i} \Big\},\nonumber\\
\forall t\in[t_{i}^{0}, t_{i}^{f}]. \label{eq:scheme1}%
\end{gather}

%%%%%%%%%%%%%%%%%%%%%%%%%%%%%%%%%%%%%%%%%%%%%%%%%%%%%%%%%%%%%%%%%%%%%%%%%%%%%%%%
%%%%%%%%%%%%%%%%%%%%%%%%%%%%%%%%%%%%
%SECTION IV: Simulation Results
%%%%%%%%%%%%%%%%%%%%%%%%%%%%%%%%%%%%

\section{Simulation Results}

The effectiveness of the efficiency of the proposed solution is validated
through simulation in VISSIM considering two intersections located in downtown
Boston. For each direction, only one lane is considered. In both intersections
the length of the control and merging zones is $L=245$ m and $S=35$ m
respectively. The distance between the two intersections is $D_{e-w}=160$ m
and $D_{w-e}=145$ m, respectively. The safe following distance is $\delta=10$
m. The arrival rate is given by a Poisson process with $\lambda=450$ veh/h.
The speed of each vehicle entering the control zone is $v_{i}^{0} = 11.11$
m/s. Note that the last two assumptions are only made for simplicity and are
by no means constraining. The desired minimum speed inside the control zones
is $7$ m/s.

\begin{figure}[ptb]
\centering
\includegraphics[width=3 in]{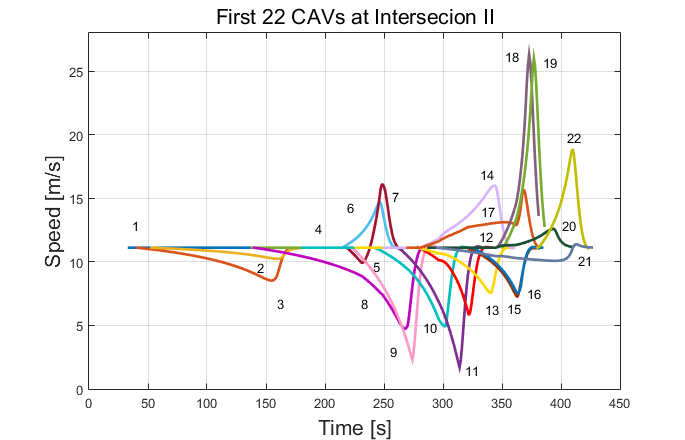} \caption{The
speed profile of the first 22 vehicles at intersection II for the
unconstrained case.}%
\label{fig:unconstrained}%
\end{figure}

\begin{figure}[ptb]
\centering
\includegraphics[width=3 in]{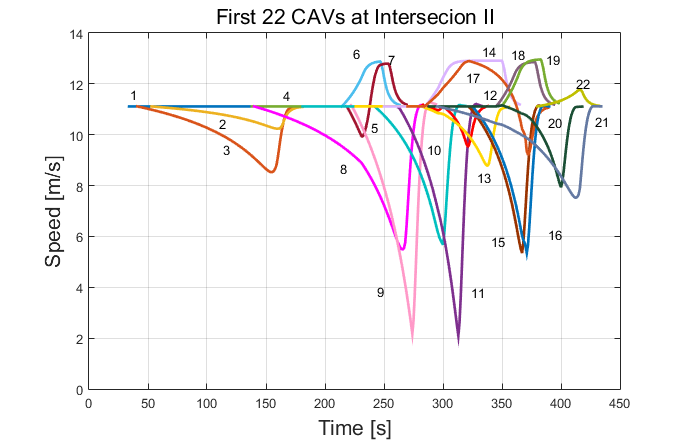} \caption{The speed
profile of the first 22 vehicles at intersection II for the constrained case.}%
\label{fig:constrained}%
\end{figure}

\begin{figure}[ptb]
\centering
\includegraphics[width=3 in]{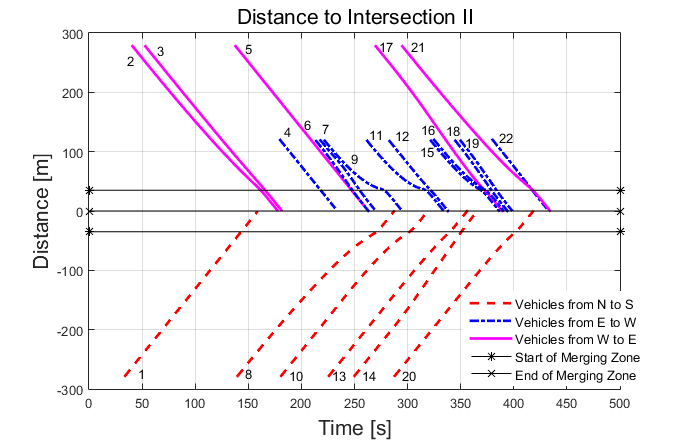}
\caption{Distance to the end of merging zone of the first 22 vehicles at
Intersection II.}%
\label{fig:position}%
\end{figure}

We considered two simulations where we: 1) relaxed the upper and lower limits
of the speed and control, and 2) included the limits. For the latter case, the
upper and lower speed limits are $v_{max}=13$ m/s and $v_{min}=0.5$ m/s
respectively. The speed of the first 22 vehicles crossing the intersection II
for both cases is shown in Fig. \ref{fig:unconstrained} and Fig.
\ref{fig:constrained} respectively. The label on each profile corresponds to
the number in the queue of the control zone II assigned by the coordinator.
The position trajectory of the first 22 vehicles crossing the intersection II
is shown in Fig. \ref{fig:position}.

Combining Fig. \ref{fig:constrained} and Fig. \ref{fig:position}, we notice
that vehicle $\#4$ assigns $t_{4}^{f} = t_{1}^{f}$ in \eqref{def:tf}. Vehicle
$\#7$ assigns $t_{7}^{f}=t_{6}^{f}+ \frac{\delta}{v_{7}(t_{6}^{f})}$ in
\eqref{def:tf}, to keep a safe distance $\delta$ from vehicle $\#6\in
\mathcal{L}_{7}$. Vehicle $\#8$ assigns $t_{8}^{f}=t_{7}^{f}+ \frac{S}%
{v_{8}(t_{7}^{f})}$ in \eqref{def:tf}, which right after vehicle $\#7
\in\mathcal{C}_{8}$ exits the merging zone, it enters the merging zone.
Vehicle $\#22$, assigns $t_{22}^{f} = t_{21}^{f}$ in \eqref{def:tf}, seems to
accelerate to catch up with vehicle $\#21$ to arrive at the merging zone at
the same time.

To investigate the interdependence between two intersections, we focus on the
behavior of three adjacent vehicles in intersection II vehicles $\#7$, $\#8$
and $\#9$, where vehicles $\#7$ and $\#9$ travel from intersection I to
intersection II. When the speed of vehicle $\#8$ becomes lower than the
desired minimum speed $7$ m/s, the vehicle $\#7$ is ``forced" to accelerate,
which creates extra space for vehicle $\#8$ to speed up (Fig.
\ref{fig:constrained} and Fig. \ref{fig:position}). Similarly, when the speed
of vehicle $\#9$ becomes lower than the minimum desired speed, vehicle $\#8$
also accelerates.

To quantify the impact of vehicle coordination on fuel consumption, we used
the polynomial metamodel \cite{Rios-Torres2015} that yields vehicle fuel
consumption as a function of the speed and acceleration. In the simulation, we
considered 448 CAVs in total crossing the two intersections and we compared
our approach to a baseline scenario that includes traffic lights with a
traffic light cycle of $30$ sec. It was shown that, with coordination of CAVs,
fuel consumption is improved by 42.4\% while the average travel time is also
improved by 37.3\% compared to the baseline scenario. The fuel consumption
improvement is due to the fact that the vehicles do not come to a full stop,
thereby conserving momentum and fuel while also improving travel time. The
average travel time for some distance inside the control zone is the same for
both cases as the vehicles approach the merging zone; however, in the baseline
scenario the vehicles need to come to a full stop resulting in increasing the
time to cross the merging zone.

%%%%%%%%%%%%%%%%%%%%%%%%%%%%%%%%%%%%%%%%%%%%%%%%%%%%%%%%%%%%%%%%%%%%%%%%%%%%%%%%
%%%%%%%%%%%%%%%%%%%%%%%%%%%%%%%%%%%%
%SECTION V: CONCLUDING REMARKS
%%%%%%%%%%%%%%%%%%%%%%%%%%%%%%%%%%%%

\section{Concluding Remarks}

The paper addressed the problem of coordinating online a continuous flow of
CAVs crossing two adjacent intersections in an urban area. We presented a
decentralized optimal control framework whose solution (when feasible) yields
for each vehicle the optimal acceleration/deceleration at any time aimed at
minimizing fuel consumption. Ongoing research investigates the feasibility of
the solution when at the time the vehicles enter the control zone some of the
constraints are active and the computational implications. Future research
should consider the diversity in CAV types crossing the intersections and the
existence of a potential trade-off between fuel consumption and congestion.

%%%%%%%%%%%%%%%%%%%%%%%%%%%%%%%%%%%%%%%%%%%%%%%%%%%%%%%%%%%%%%%%%%%%%%%%%%%%%%%%
\bibliographystyle{IEEETran}
\bibliography{TITS_merging}

\begin{thebibliography}{10}
\providecommand{\url}[1]{#1}
\csname url@rmstyle\endcsname
\providecommand{\newblock}{\relax}
\providecommand{\bibinfo}[2]{#2}
\providecommand\BIBentrySTDinterwordspacing{\spaceskip=0pt\relax}
\providecommand\BIBentryALTinterwordstretchfactor{4}
\providecommand\BIBentryALTinterwordspacing{\spaceskip=\fontdimen2\font plus
\BIBentryALTinterwordstretchfactor\fontdimen3\font minus
  \fontdimen4\font\relax}
\providecommand\BIBforeignlanguage[2]{{%
\expandafter\ifx\csname l@#1\endcsname\relax
\typeout{** WARNING: IEEEtran.bst: No hyphenation pattern has been}%
\typeout{** loaded for the language `#1'. Using the pattern for}%
\typeout{** the default language instead.}%
\else
\language=\csname l@#1\endcsname
\fi
#2}}

\bibitem{Li2014}
L.~Li, D.~Wen, and D.~Yao, ``{A Survey of Traffic Control With Vehicular
  Communications},'' \emph{IEEE Transactions on Intelligent Transportation
  Systems}, vol.~15, no.~1, pp. 425--432, 2014.

\bibitem{Liu2007}
Z.~Liu, ``A survey of intelligence methods in urban traffic signal control,''
  \emph{IJCSNS International Journal of Computer Science and Network Security},
  vol.~7, no.~7, pp. 105--112, 2007.

\bibitem{Fleck2015}
J.~L. Fleck, C.~G. Cassandras, and Y.~Geng, ``Adaptive quasi-dynamic traffic
  light control,'' \emph{IEEE Transactions on Control Systems Technology},
  2015, DOI: 10.1109/TCST.2015.2468181, to appear.

\bibitem{Dresner2004}
K.~Dresner and P.~Stone, ``{Multiagent traffic management: a reservation-based
  intersection control mechanism},'' in \emph{Proceedings of the Third
  International Joint Conference on Autonomous Agents and Multiagents Systems},
  2004, pp. 530--537.

\bibitem{Dresner2008}
------, ``{A Multiagent Approach to Autonomous Intersection Management},''
  \emph{Journal of Artificial Intelligence Research}, vol.~31, pp. 591--653,
  2008.

\bibitem{DeLaFortelle2010}
A.~{de La Fortelle}, ``{Analysis of reservation algorithms for cooperative
  planning at intersections},'' \emph{13th International IEEE Conference on
  Intelligent Transportation Systems}, pp. 445--449, Sept. 2010.

\bibitem{Huang2012}
S.~Huang, A.~Sadek, and Y.~Zhao, ``{Assessing the Mobility and Environmental
  Benefits of Reservation-Based Intelligent Intersections Using an Integrated
  Simulator},'' \emph{IEEE Transactions on Intelligent Transportation Systems},
  vol.~13, no.~3, pp. 1201,1214, 2012.

\bibitem{Milanes2010}
V.~Milan\'{e}s, J.~P\'{e}rez, and E.~Onieva, ``{Controller for Urban
  Intersections Based on Wireless Communications and Fuzzy Logic},'' \emph{IEEE
  Transactions on Intelligent Transportation Systems}, vol.~11, no.~1, pp.
  243--248, 2010.

\bibitem{Alonso2011}
J.~Alonso, V.~Milan\'{e}s, J.~P\'{e}rez, E.~Onieva, C.~Gonz\'{a}lez, and
  T.~de~Pedro, ``{Autonomous vehicle control systems for safe crossroads},''
  \emph{Transportation Research Part C: Emerging Technologies}, vol.~19, no.~6,
  pp. 1095--1110, Dec. 2011.

\bibitem{Kim2014}
K.-D. Kim and P.~Kumar, ``{An MPC-Based Approach to Provable System-Wide Safety
  and Liveness of Autonomous Ground Traffic},'' \emph{IEEE Transactions on
  Automatic Control}, vol.~59, no.~12, pp. 3341--3356, 2014.

\bibitem{Campos2014}
G.~R. Campos, P.~Falcone, H.~Wymeersch, R.~Hult, and J.~Sjoberg, ``Cooperative
  receding horizon conflict resolution at traffic intersections,'' in
  \emph{2014 IEEE 53rd Annual Conference on Decision and Control (CDC)}, 2014,
  pp. 2932--2937.

\bibitem{Malikopoulos2016}
A.~A. Malikopoulos and C.~G. Cassandras, ``Decentralized optimal control for
  connected and automated vehicles at an intersection,'' in \emph{55th
  Conference on Decision and Control}, 2016, - arXiv:1479353.

\bibitem{ZhangMalikopoulosCassandras2016}
Y.~Zhang, A.~A. Malikopoulos, and C.~G. Cassandras, ``Optimal control and
  coordination of connected and automated vehicles at urban traffic
  intersections,'' 2016,- arXiv:1362458.

\bibitem{Rios-Torres2015}
J.~Rios-Torres, A.~A. Malikopoulos, and P.~Pisu, ``{Online Optimal Control of
  Connected Vehicles for Efficient Traffic Flow at Merging Roads},'' in
  \emph{2015 IEEE 18th International Conference on Intelligent Transportation
  Systems}, Canary Islands, Spain, September 15-18, 2015.

\bibitem{Zhong2010}
M.~Zhong and C.~G. Cassandras, ``Asynchronous distributed optimization with
  event-driven communication,'' \emph{IEEE Transactions on Automatic Control},
  vol.~55, no.~12, pp. 2735--2750, 2010.

\end{thebibliography}

\end{document}